\documentclass[11pt]{article}
\usepackage{amssymb, amsmath, amsthm, textcomp}

\theoremstyle{plain}
\begin{document}
\parindent 0in
\parskip 1 em
\date{}
\title{\text{\bf{\normalsize{TORSION-FREE $G_2$-STRUCTURES WITH IDENTICAL RIEMANNIAN METRIC}}}}
\author{Christopher Lin}
\maketitle 

\begin{abstract}
 Based on a general formula due to R.Bryant, we work out the topological structure of the space of 
 torsion-free $G_2$-structures generating the same associated Riemannian metric on a compact $7$-manifold.  We also identify a corresponding Lie group-theoretic structure of the space.  These observations are then used to describe the moduli space of torsion-free $G_2$-structures in certain cases - by way of covering spaces.        
\end{abstract}

Preprint of an article accepted for publication in \text{\it{Journal of Topology and Analysis}} (2018),\\ 
DOI: 10.1142/S1793525318500310 \textcircled{c} World Scientific Publishing Company, \\
www.worldscientific.com/worldscinet/jta.  

\section{Introduction}
The group $G_2$ can be defined as the group of automorphisms of the Octonions $\mathbb{O}_8$.  This definition can be shown to be equivalent 
to the subgroup of $SO(7)$ that preserves the $3$-form 
\begin{equation}\label{phi0}
 \varphi_0 = dx_{123}+dx_{145} + dx_{167}+dx_{246} -dx_{257} - dx_{347}-dx_{356}
\end{equation}  
in $\mathbb{R}^7$, where $dx_{ijk} = dx_i \wedge dx_j \wedge dx_k$.  On a $7$-dimensional manifold $M$, a $G_2$-structure is a principal sub-bundle 
of the $GL(7,\mathbb{R})$-frame bundle over $M$, with the reduced structure group $G_2$.  Analytically, a $G_2$-structure is described by a 
smooth $3$-form $\varphi$ on $M$ resembling (\ref{phi0}) in appropriate local frames.  The switching of thinking between the two equivalent ways of 
describing $G_2$-structures is a theme throughout this paper.    

\quad Because $G_2$ is a subgroup of $SO(7)$, a $G_2$-structure is contained in a unique $SO(7)$-structure, i.e. it is associated to 
a unique Riemannian metric together with a unique orientation.  In the analytic language, we say that a $G_2$-structure $\varphi$ \text{\it{generates}} 
a unique Riemannian metric $g_\varphi$.    
However, different $G_2$-structures may generate the same 
metric.  It is therefore natural to ask: for a fixed Riemannian metric generated by some $G_2$-structure, 
what are the distinct $G_2$-structures that generate the same metric?  Equivalently, 
we are asking for a description of the space of $G_2$-structures contained in a single $SO(7)$-structure over a $7$-manifold.  
R. Bryant gave an answer to this general question in 
\cite{Br}: given a $G_2$-structure 
$\varphi$ on a $7$-manifold $M$, for any $f\in C^{\infty}(M)$ and 
$\omega\in\Omega^1(M)$ such that $f^2 + |\omega|^2 = 1$ the $3$-form $\tilde{\varphi}$ given by 
\begin{equation}\label{claim}
 \tilde{\varphi} = (f^2 - |\omega|^2)\varphi + 2f\ast(\omega\wedge\varphi) + 2\omega\wedge\ast(\omega\wedge\ast\varphi),
\end{equation} 
is a $G_2$-structure that generates the same Riemannian metric and orientation as $\varphi$, where the Hodge star $\ast$ and norm 
$|\cdot|$ are induced by the metric $g_\varphi$.  Moreover, any $G_2$-structure $\tilde{\varphi}$ that induces the same metric and orientation as 
$\varphi$ can be written as (\ref{claim}) for some 
$f\in C^{\infty}(M)$ and $\omega\in\Omega^1(M)$ such that $f^2 + |\omega|^2 = 1$, unique up to replacement by $(-f,-\omega)$.  Therefore, the space of $G_2$-structures that generate the same Riemannian metric is equivalent to the space of sections of a $\mathbb{RP}^7$-bundle over $M$.

\quad On the other hand, a \text{\it{torsion-free}}
$G_2$-structure is one which is preserved by parallel translation with respect to the metric it generates.  Therefore there should be a more elegant description of the finer space 
of all torsion-free $G_2$-structures 
that generate the same metric and orientation.  To this end, let $\varphi$ be a torsion-free $G_2$-structure on $M$ and we define the space
\begin{equation*}
 \Gamma = \{\text{torsion-free $G_2$-structures on $M$ that generate $g_\varphi$, with the same orientation as $\varphi$}\}.
\end{equation*}
Canonically we think of $\Gamma$ as a set of $3$-forms on $M$.  It turns out that the topological structure of $\Gamma$ is 
independent of $\varphi$.  Since this paper mainly concerns the topology of $\Gamma$, we omit $\varphi$ in the notation for $\Gamma$.

\quad The main goal of this paper is to derive an analytic description of 
$\Gamma$ from (\ref{claim}) that seems to be missing in the literature.  For this part of the paper we shall assume that $M$ is compact.  
In particular, we will see 
that this description implies $\Gamma$ is topologically $\mathbb{RP}^{\,b^1}$, where 
$b^1$ is the first Betti number of $M$.  In fact we will show that $\Gamma$ is $C^1$-diffeomorphic to $\mathbb{RP}^{b^1}$ (see Proposition \ref{except}).  
It turns out that the topology of $\Gamma$, along with a particular group of isometries $S_0$ to be defined later, determines the topology in a particular direction on the moduli space $\mathcal{M}$ of torsion-free $G_2$-structures on a compact $7$-manifold.  Denote by 
$X$ the projected image of $\Gamma$ into $\mathcal{M}$, we will prove the following result. 

\newtheorem{theorem}{Theorem}[section]
\begin{theorem}\label{Thm1}
 The projection $\Gamma \longrightarrow X$ is a covering map.
\end{theorem}
Various consequences of Theorem \ref{Thm1} on the moduli space $\mathcal{M}$ will be discussed for special cases.  As $\Gamma$ is homeomorphic to real projective space, 
understanding the topology of $X$ is inevitably related to the classic problem of classifying free group actions on spheres.  If $b^1 >1$, we will see that covering space theory also implies that, provided $\Gamma\longrightarrow X$ is a normal covering, the fundamental group of $X$ must be a central extension of the group of deck transformations of $\Gamma\longrightarrow X$ by $\mathbb{Z}_2$ - the fundamental group of real projective space of dimension bigger than $1$.  A complete understanding of the topoology of $X$ hinges on a full understanding of these problems in algebra/algebraic-topology.

\quad We would like to point out that the global topology of $\mathcal{M}$ is still unknown, and is a very important problem.  In particular, geometers and physicists are very interested in the case when the underlying $M$ has full $G_2$-holonomy.  In this case, $b^1=0$ and $X$ is just a point.  Thus our results are only interesting when the underlying manifold has holonomy group properly contained in $G_2$.  However, because of the availability of concrete examples when the holonomy group is properly contained in $G_2$, we believe our results still have merit in understanding the subject as a whole.              
        
\quad On the other hand, there is a general Lie group-theoretic description of $\Gamma$ which does not require $M$ to be compact.  To be more precise, we choose 
a point $p\in M$ over which a frame $f$ is adapted to a torsion-free $G_2$-structure in $\Gamma$.  
Then we define the set
\begin{equation}\label{N_f}
 N_f = \{g\in SO(7)\,|\,g^{-1}\text{Hol}_f(D)g \subset G_2\},
\end{equation}  
where $\text{Hol}_f(D)$ denotes the holonomy group (based at $f$) of the Levi-Civita connection $D$ associated to the unique Riemannian metric of $\Gamma$.  For any chosen point-wise frame $f$, denote by $[f]$ the orbit of the right-action by elements of $G_2$.  Then $[f]$ represents a torsion-free $G_2$-structure (hence an element of $\Gamma$) by parallel translation with respect to the induced Riemannian metric.  Then every element in $\Gamma$ can be denoted by $[f.h]$ for some $h\in SO(7)$, where $f.h$ denotes the 
canonical right-action of $h$ on the frame $f$.  We shall demonstrate explicitly in Section \ref{Lie} that 
\newtheorem{proposition}[theorem]{Proposition}  
 \begin{proposition}\label{prop1}
  The map $N_f / G_2 \longrightarrow \Gamma$ given by $h\,G_2 \mapsto [f.h]$ is a bijection.  
 \end{proposition}
In other words, $\Gamma$ is parameterized by the coset space $N_f/G_2$.  It should be pointed out that Proposition \ref{prop1} is actually a special case of Proposition 3.1.8 in \cite{J1}, although it is often ignored.  Via our first result, we will show that $N_f/G_2$ is in fact also homeomorphic to 
$\mathbb{RP}^{\,b^1}$ when $M$ is compact.  An important message here is that the respective topologies of $\Gamma$ and $N_f / G_2$ only depend on the topology of the 
underlying manifold $M$.     

\quad The organization of this paper is as follows.  In the next section, we will discuss the necessary  background concerning $G_2$-structures on $7$-manifolds and the associated moduli space.  In Section \ref{Bry}, we will derive the analytic description of $\Gamma$ from Bryant's formula (\ref{claim}) and identify the topological structure of $\Gamma$.  
In Section \ref{Lie}, we will clarify the characterization of $\Gamma$ that uses frames, prove the homeomorphism of 
$\Gamma$ and $N_f / G_2$ and discuss its relation to the analytic description.  The background material in Section \ref{prelim} is presented from the analytic point of view, while the equivalent principal bundle version will appear in Section \ref{Lie} since it is more relevant there.  Finally, in the last section 
we will prove Theorem \ref{Thm1} and discuss the insights it provides into the topology of the moduli space 
$\mathcal{M}$ in special cases.  The associated free action on spheres and group extension problems will highlight the end of the paper.                


\section{Preliminaries on $G_2$-structures}\label{prelim}
A $G_2$-structure on a $7$-manifold $M$ is a $3$-form $\varphi$ such that there are local tangent frames $f$ for which $\varphi = f^{*}\varphi_0$, with $\varphi_0$ 
as in (\ref{phi0}).  Throughout the paper $M$ is always connected.  For each such local frame $f$ and a point $p\in M$ in the relevant trivializing neighborhood, denote by $f_p$ as the value of the frame at $p$.  Here we are treating $f_p=\{e_i\}_{i=1}^7$ as the point-wise linear map $f:T_pM\longrightarrow \mathbb{R}^7$ given by $f(e_i) = \vec{e}_i$, where $\{e_i\}_{i=1}^7$ is a basis $T_pM$ and $\vec{e}_i$ is the natural unit vector.  The transition functions of these local frames lie exactly in the group $G_2$ as defined via (\ref{phi0}).  $G_2$ acts naturally 
on $\mathbb{R}^7$ and its exterior algebra, giving rise to irreducible representations.  Given a $G_2$-structure $\varphi$ on $M$, these irreducible representations pass onto the spaces of differential forms.  In particular, we have 
\begin{align}
 \Omega^2(M) &= \Omega^2_7 \oplus \Omega^2_{14} \notag\\
 \Omega^3(M) &= \Omega_1^3 \oplus \Omega_7^3 \oplus \Omega_{27}^3 \notag
\end{align}
where 
\begin{align}\label{decomp}
 \Omega^2_7 &= \{X\lrcorner\varphi \,|\, X\in \Gamma(TM)\} =\{\beta\in\Omega^2(M)\,|\,\ast(\varphi\wedge\beta)=-2\beta\} \notag\\
 \Omega^2_{14} &= \{\beta\in\Omega^2(M)\,|\, \beta\wedge\ast\varphi = 0\}\notag\\
 \Omega^3_1 &= \{f\varphi\,|\,f\in C^{\infty}(M)\} \notag\\
 \Omega^3_7 &= \{X\lrcorner\ast\varphi\,|\,X\in\Gamma(TM)\}\notag\\
 \Omega^3_{27} &= \{\eta\odot\varphi\,|\,\eta \,\,\text{trace-less symmetric 2-tensor}\},
\end{align}
where $\odot$ is defined by the canonical linear map $T^{*}M\otimes T^{*}M\longrightarrow \Omega^3(M)$ given by   
\begin{equation}\label{odot}
(b\odot\varphi) (u,v,w) = \varphi(\hat{b}u,v,w) + \varphi(u,\hat{b}v,w) + \varphi(u,v,\hat{b}w)
\end{equation} with $\hat{b}\in \,\text{End}(TM)$ the $g_\varphi$-dual of 
$b\in T^{*}M\otimes T^{*}M$.  The corresponding decompositions for $\Omega^4(M)$ and $\Omega^5(M)$ are defined by taking the Hodge star of the irreducible components of 
$\Omega^3(M)$ and $\Omega^2(M)$, respectively.  More generally, the map (\ref{odot}) 
encodes a fundamental relation between $\Omega^2(M)$ 
and $\Omega^3(M)$, as seen by the decomposition
\begin{align}\label{split}
 T^{*}M\otimes T^{*}M &= \text{Sym}^2(TM)\oplus\Omega^2(M) \notag\\
                     &= \Big(C^{\infty}(M)\otimes g_\varphi\Big) \oplus \text{Sym}^2_0(TM) \oplus\Omega^2_7\oplus\Omega^2_{14} 
                     \longrightarrow \Omega_1^3 \oplus \Omega^3_7 \oplus \Omega^3_{27},
\end{align}
where $C^{\infty}(M)\otimes g_\varphi$ is mapped isomorphically onto $\Omega_1^3$, the space of trace-less symmetric $2$-tensors $\text{Sym}^2_0(TM)$ 
is mapped isomorphically onto $\Omega^3_{27}$, 
$\Omega^2_7$ is mapped isomorphically onto $\Omega^3_7$, and $\Omega^2_{14}$ is the kernel of the map.  
In fact, $\Omega^2_{14}$ is isomorphic to the Lie Algebra of $G_2$.  We would also like to point out that with respect to any (point-wise) frame $\{e_i\}_{i=1}^7$ adapted to the $G_2$-structure, it can be readily verified from (\ref{odot}) that 
\begin{align}\label{Jazz}
 b\odot\varphi &= b_{ij}\,\omega_i\wedge(e_j\lrcorner\varphi) \notag\\
 &= b_{ij}\,\omega_i\wedge\ast(\omega_j\wedge\ast\varphi),
\end{align}   
where $\omega_i$ is the metric dual 1-form of $e_i$ 
and $b_{ij} = b(e_i,e_j)$.  This local expression is useful in many situations.  In fact, we will make use of it in Section \ref{Bry}.

\quad The necessary and sufficient conditions for the existence of a $G_2$-structure on a $7$-manifold is that it be orientable and spin, which are 
relatively mild topological conditions.  What we are interested in are \text{\it{torsion-free}} $G_2$-structures, which are a lot harder to 
come by.  A $G_2$-structure $\varphi$ is said to be torsion-free if $\nabla\varphi = 0$ on $M$, where $\nabla$ is the Levi-Civita connection 
corresponding to $g_\varphi$.  It follows immediately that $d\varphi =0 = \delta\varphi$ if $\varphi$ is torsion-free.  It was the work of 
Fernandez and Gray \cite{FG} that $d\varphi =0 = \delta\varphi$ also implies $\nabla\varphi=0$.  If $M$ is compact $d\varphi =0 = \delta\varphi$ 
is also equivalent to $\Delta\varphi=0$, where $\Delta = d\delta+\delta d$ is the Hodge Laplacian induced by $g_\varphi$.  The torsion-free condition on $\varphi$ is also equivalent to the Riemannian holonomy group of $g_\varphi$ being contained in $G_2$.  We say $M$ is a 
$G_2$-manifold if it admits a torsion-free $G_2$-structure.  It can be shown that a $G_2$-manifold is always Ricci-flat.  

\quad The most important consequence of having a torsion-free $G_2$-structure is that the decomposition (\ref{decomp}) descends to cohomology when 
$M$ is compact.  Throughout this paper, all cohomology groups have coefficients in $\mathbb{R}$.  If 
$\varphi$ is torsion-free, then the corresponding Hodge Laplacian $\Delta$ maps $\Omega^k_l$ to itself.  Therefore via the Hodge Theorem we can define 
the cohomology group $H^k_l(M)$, which is a subspace of $H^k(M)$ whose elements are represented by harmonic $\Omega^k_l$-forms.  We then have
\begin{align}\label{coho}
 H^2(M) &= H^2_7(M) \oplus H^2_{14}(M) \notag\\
 H^3(M) &= H^3_1(M) \oplus H^3_7(M) \oplus H^3_{27}(M)
\end{align}    
with the other cohomology groups defined via the Hodge star operator.  It can be shown that $H^1(M)\cong H^k_7(M)$ for $k=2,3,4,5$.  Defining 
the \text{\it{reduced Betti numbers}}\, $b_l^k = \text{dim}\, H^k_l(M)$, we thus have $b_7^k = b^1$ and that the only unknown ones are 
essentially $b_{14}^{2}$ and $b_{27}^{3}$.  Note that $b_1^3 = 1$ always.  It follows that $b_{14}^{2}$ and $b_{27}^3$ are topological as well.  We also know that for a compact $G_2$-manifold, 
$\text{Hol}_f(D) = G_2$ if and only if the fundamental group of $M$ is finite (see \cite{J1}), and therefore $\text{Hol}_f(D) = G_2$ implies $b^1=0$. 



\quad One of the most seminal results in $G_2$-geometry was D. Joyce's construction of compact $G_2$-manifolds of full $G_2$-holonomy (see his book \cite{J1}).  The examples he painstakingly constructed were non-explicit, due to the fact that such manifolds do not have any 
``$G_2$-symmetry'' (see \cite{Lin} for an explanation).  On the other hand, there are well-known concrete examples when the holonomy group is properly contained in $G_2$:

 \begin{enumerate}
  \item $S^1 \times Y$ with $\varphi = 
  d\theta\wedge\omega+\text{Re}\,\Omega$, where $Y$ is a 
  Calabi-Yau $3$-fold with K\"ahler form $\omega$ and 
  $\Omega$ its holomorphic volume form.  The holonomy group is $SU(3)$.  
  \item $\mathbb{T}^3\times Y$ with $\varphi = dx_{123}+
  dx_1\wedge\omega +dx_2\wedge\text{Re}\,\Omega - 
  dx_3\wedge\text{Im}\,\Omega$, where $Y$ is a Calabi-Yau 
  $2$-fold with K\"ahler form $\omega$ and $\Omega$ its holomorphic volume form.  The holonomy group is $SU(2)$.
  \item $\mathbb{T}^7$ with the inherited $G_2$-structure (\ref{phi0}) from $\mathbb{R}^7$.  The holonomy group is 
  $\{1\}$.
 \end{enumerate}
 
As shown in \cite{Lin}, the dimension of $G_2$-symmetries of a torsion-free $G_2$-structure on a compact manifold is equal to the first Betti number.  For the three examples above they are $1, 3, 7$ in order, and such symmetries are exhibited explicitly by the respective torus parts of the manifold.  

\quad Let us end this section by recalling the moduli space of torsion-free $G_2$-structures.  
In \cite{J1}, D. Joyce defined the Moduli space $\mathcal{M} = \{\varphi\,\,\text{torsion-free of same orientation}\} / \,\text{Diff}_0(M)$ 
on a compact $7$-manifold $M$, where $\text{Diff}_0(M)$ denotes the group of diffeomorphisms of $M$ to itself that are 
isotopic to the identity.  Let $\{\varphi\}\in\mathcal{M}$ denote the orbit of the action of $\text{Diff}_0(M)$ on $\varphi$ by pull-back.  
Note that the action by diffeomorphisms isotopic to the identity ensures that all elements in $\{\varphi\}$ belong to the same cohomology class 
$[\varphi]$ in $H^3(M)$.  It is proved in \cite{J1} that $\mathcal{M}$ is locally diffeomorphic to $H^3(M)$ through the map $\{\varphi\}\mapsto [\varphi]$.    
However, the global topology of $\mathcal{M}$ is unknown.  On the other hand, $c\varphi$ is still a torsion-free 
$G_2$-structure with the same orientation as $\varphi$, for any $c>0$.  Therefore any local chart given 
by Joyce's map above can be extended towards $0\in H^3(M)$ and away to infinity along the $H^3_1(M)$-direction tangential to $\mathcal{M}$.  
Hence $\mathcal{M}$ actually has a cone structure - with the unknown topology being in directions tangent to the 
subspace $H^3_7(M)\oplus H^3_{27}(M)$ 
(the ``collars'' of the cone) at each point on $\mathcal{M}$.  What is also unknown is the number of connected components of $\mathcal{M}$.  We will eventually see in Section \ref{Moduli} 
that the topology in the $H^3_7(M)$-direction at every point on $\mathcal{M}$ is manifested as the orbit space of a group consisting of isometries acting on $\Gamma$.    

\section{Analytic description based on Bryant's formula}\label{Bry}
In this section, we will derive the following analytic description of $\Gamma$ from formula (\ref{claim}). 
\newtheorem{prop3}{Proposition}[section]
\begin{prop3}\label{Bryant}
On a compact $7$-manifold endowed with a torsion-free $G_2$-structure $\varphi$, there is the bijection
 \[
 \Gamma \cong \Big\{(c,\omega)\in [-1,1]\times H^1(M) \,\,\Big|\,\, c^2 + |\omega|^2 = 1\Big\} \Big/ \Big\{(c,\omega)\sim (-c,-\omega)\Big\},
 \]
 where $\omega$ denotes the harmonic representative of an element of $H^1(M)$ with respect to $g_\varphi$.
\end{prop3}  

Note that since $\omega$ is a harmonic $1$-form on a compact, Ricci-flat manifold, it is also parallel. Therefore the point-wise norm $|\omega|$ is constant on $M$, and hence the condition 
$c^2 + |\omega|^2 = 1$ makes sense on $M$.  In particular, up to a constant multiple $|\omega|$ serves as a norm on $H^1(M)$.      

\quad Before giving the proof of Proposition \ref{Bryant}, let us provide some extra motivation.  Let $\varphi(t)$ be a smooth family of $G_2$-structures.  Then according to 
(\ref{decomp}), (\ref{odot}), and (\ref{split}) we can write 
\[
 \frac{\partial \varphi}{\partial t} = b(t)\odot\varphi(t) + X(t)\lrcorner\ast\varphi(t)
\]
for a unique family of symmetric $2$-tensors $b(t)$ and a unique family of vector fields $X(t)$, and where $\odot$ and $\ast$ 
are induced by $\varphi(t)$.  It can be shown (see \cite{K1}) that 
\[
 \frac{\partial g}{\partial t} = 2b(t),
\]  
where $g(t) = g_{\varphi(t)}$ is the  corresponding family of metrics induced.  Now suppose $\varphi(t)$ generate the same metric for all $t$, then it follows 
that $\frac{\partial\varphi}{\partial t} \in \Omega^3_7$ for all $t$.  If in addition $\varphi(t)$ is torsion-free for all $t$, we see that
\[
  d\frac{\partial\varphi(t)}{\partial t} = \frac{\partial}{\partial t}d\varphi(t) = 0
\] 
and 
\[
  d\ast\frac{\partial\varphi(t)}{\partial t} = \frac{\partial}{\partial t}d\ast\varphi(t) = 0,
\]   
where the fixed metric also fixes the Hodge star $\ast$.  It follows that $\frac{\partial\varphi}{\partial t}$ is harmonic for all $t$ and 
thus represents an element in $H^3_7(M)$.  On the moduli space $\mathcal{M}$ this means that variations which preserve the metric are always in 
the $H^3_7(M)$-directions.  It follows that $\Gamma$ is 
a discrete set in $\mathcal{M}$ if $H^1(M)=0$ (which happens when $M$ has full $G_2$-holonomy).  
This local result on $\mathcal{M}$ is refined by Proposition \ref{Bryant}, which says that in fact $\Gamma$ consists of a single point.  In the next section, we will see yet another way to arrive at this fact.        

\quad On the other hand, suppose $H^3_7(M)\ne 0$.  At any point on $\mathcal{M}$ and given a nonzero element $[\eta]$ of $H^3_7(M)$, 
is there always a variation of torsion-free $G_2$-structures preserving 
the metric in the direction of $[\eta]$?  After the proof of Proposition \ref{Bryant} below, we will prove that the bijection there is actually a $C^1$-diffeomorphism, which consequently means that $\Gamma$ is homeomorphic to $\mathbb{RP}^{\,b^1}$.  
In the last section we will show that $\Gamma$ projects into $\mathcal{M}$ via a covering map.  Thus locally, the image of $\Gamma$ in 
$\mathcal{M}$ still has dimension $b_1$.  Then since $b^3_7 = b^1$, this answers the question above in the affirmative.

\text{\bf{Proof of Proposition \ref{Bryant}.}} \quad First note that the first two terms on the right-hand side of (\ref{claim}) are respectively in $\Omega^3_1$ and 
$\Omega^3_{7}$.  For the third term however, we note that 
\begin{align}
 \omega\wedge\ast(\omega\wedge\ast\varphi)\wedge\ast(X\lrcorner\ast\varphi) 
 &= -\omega\wedge\ast(\omega\wedge\ast\varphi)\wedge X^{\flat}\wedge\varphi \notag\\
 &= -\omega\wedge X^{\flat}\wedge\varphi\wedge\ast(\omega\wedge\ast\varphi) \notag\\
 &= 2\omega\wedge X^{\flat}\wedge\omega\wedge\ast\varphi \notag\\
 &=0,
\end{align}
where in arriving at the third equality we observed that $\ast(\omega\wedge\ast\varphi)\in\Omega^2_7$ 
(see (\ref{decomp})). It follows that the third term has no $\Omega^3_7$-component.  On the other hand, let 
$\pi_{1}\big(\omega\wedge\ast(\omega\wedge\ast\varphi)\big) = \alpha\, \varphi$ for some function $\alpha$ on $M$, and we see that 
\begin{align}
 7\alpha\,dv_{\varphi} &= \alpha \,\varphi\wedge\ast\varphi \notag\\ 
                       &= \omega\wedge\ast(\omega\wedge\ast\varphi)\wedge\ast\varphi \notag\\
                       &= \omega\wedge\ast 3\omega \notag\\
                       &= 3|\omega|^2 dv_{\varphi}, \notag
\end{align}
where in arriving at the third equality we used the standard identity $\ast\varphi\wedge\ast(\ast\varphi\wedge\alpha) = 3\ast\alpha$ for any $\alpha\in\Omega^1(M)$ (see Proposition A.3 of \cite{K1}).
Thus $\alpha = \frac{3}{7}|\omega|^2$, which is in general nonzero.  
Hence $\omega\wedge\ast(\omega\wedge\ast\varphi)\in\,\Omega^3_1\oplus\Omega^3_{27}$.  We can then rewrite (\ref{claim}) explicitly in its 
$\Omega^3_k$-components as 
\begin{align}\label{reborn}
 \tilde{\varphi} &= (f^2 - |\omega|^2)\varphi + 2f\ast(\omega\wedge\varphi) + 2\Big(\frac{3}{7}|\omega|^2\varphi + \pi_{27}\big(\omega\wedge\ast(\omega\wedge\ast\varphi)\big)\Big) \notag\\
 &= \frac{1}{7}(8f^2-1)\varphi +  2f\ast(\omega\wedge\varphi)+ 2\pi_{27}\big(\omega\wedge\ast(\omega\wedge\ast\varphi)\big),
\end{align}  
where we have also used $f^2+|\omega|^2 = 1$.   

\quad Now, $\tilde{\varphi}$ is torsion-free, and therefore harmonic with respect to $g_{\tilde{\varphi}} = g_\varphi$.  Taking 
the Laplacian with respect to $g_\varphi$ on 
(\ref{reborn}) and using the fact that the Laplacian commutes with projections to irreducible subspaces, it follows that $8f^2 - 1$ must be harmonic and hence $f$ must be constant on the compact manifold $M$. 

\quad On the other hand, because $f=c$ is constant and $\varphi$ is torsion-free, we also see that
\begin{align}
 0 &= \Delta\ast(f\omega\wedge\varphi) \notag\\
   &= \ast c\,\Delta(\omega\wedge\varphi) \notag\\
   &= c\ast\big(\Delta\omega\wedge\varphi\big). \notag
\end{align} 
Thus $\Delta\omega =0$ necessarily. 

\quad Next, let $(c,\omega)\in [-1,1]\times \mathcal{H}^1(M)$.  Then by taking the Laplacian again, it is clear 
that the first two terms of the right-hand side of (\ref{claim}) go to zero.  For the last term, 
since $\Delta\omega = 0$ on $M$ and $M$ is compact, it follows again by the Ricci-flatness of $M$ that $\nabla\omega=0$ on $M$.  Since $\nabla\varphi=0$ on $M$ as well, by the product rule and the fact that $\nabla$ commutes with $\ast$, it follows 
that $\omega\wedge\ast(\omega\wedge\ast\varphi)$ is parallel as well.  This then implies that 
$\omega\wedge\ast(\omega\wedge\ast\varphi)$ must be both closed and co-closed, hence harmonic.  Therefore $(c,\omega)$ indeed represents an element in $\Gamma$.

\quad Finally, suppose $(c_1 , \omega_1),\, (c_2, \omega_2)\in  [-1,1]\times \mathcal{H}^1(M)$ both determine the same $\tilde{\varphi}\in \Gamma$.  
Note from (\ref{reborn}) that $\pi_1\big(\tilde{\varphi}\big) = \frac{1}{7}(8c^2 - 1)\varphi$ and $\pi_7\big(\tilde{\varphi}\big)=2c\ast(\omega\wedge\varphi)$ for some $c\in [-1,1]$ and 
$\omega\in\mathcal{H}^1(M)$.  By comparing the corresponding irreducible components we see that $c_1 =\pm c_2$ and $\omega_1 = \pm\omega_2$, 
therefore we indeed have a bijection between $\Gamma$ and the set in question.  
  \qed
  
Next we investigate the regularity of the map above.

\begin{prop3}\label{except}
 The map in Proposition \ref{Bryant} is a $C^1$-diffeomorphism, hence a homeomorphism.  
\end{prop3}

\text{\bf{Proof.}}\quad  The map $(c,\omega)\mapsto \tilde{\varphi}(c,\omega)$ in Proposition 
\ref{Bryant} given by 
\begin{equation}\label{aroma}
 \tilde{\varphi}(c,\omega)= 
 (c^2 - |\omega|^2)\varphi + 
 2c\ast(\omega\wedge\varphi) + 
 2\omega\wedge\ast(\omega\wedge\ast\varphi), 
\end{equation}
was shown to be bijective onto $\Gamma$ up to the equivalence $(c,\omega) \sim (-c,-\omega)$, and is clearly $C^1$ as a map into the finite-dimensional space of harmonic $3$-forms with respect to $g_\varphi$.  We will show that the induced derivative map from the tangent space of any point on the $b^1$-dimensional sphere \[
\Big\{(c,\omega)\in [-1,1]\times H^1(M) \,\,\Big|\,\, c^2 + |\omega|^2 = 1\Big\}
\] is injective.  
Then by the inverse function theorem, this shows that the map (\ref{aroma}) is a $C^1$-diffeomorphism onto $\Gamma$ up to the 
equivalence $(c,\omega) \sim (-c,-\omega)$ with a $C^1$ inverse.  Ultimately the map is a topological embedding and hence a homeomorphism onto $\Gamma$.  

\quad Let $(c_0,\omega_0)$ be a point on the sphere above, and let $(c(t), \omega(t))$ be a $C^1$-curve on the sphere such that $(c(0), \omega(0))= (c_0, \omega_0)$.  In view of 
formula (\ref{aroma}), suppose
\begin{align}\label{metaji}
 \frac{\partial}{\partial t}
 \tilde{\varphi}(c(t),\omega(t))\Big|_{t=0} &= 
 4c_0\,\dot{c}(0)\varphi \,+ \, 2\dot{c}(0)\ast(\omega_0\wedge\varphi)\,+\, 2 c_0\ast(\dot{\omega}(0)\wedge\varphi) \notag\\
 &\hskip 0.5cm+ 
 2\dot{\omega}(0)\wedge\ast(\omega_0\wedge\ast\varphi) 
 \, +\, 2\omega_0\wedge\ast(\dot{\omega}(0)\wedge\ast\varphi)\notag\\
 &=0. 
\end{align}
Note that we have used the condition $c^2 + |\omega|^2=1$ in deriving the term $4c_0\,\dot{c}(0)$ above.  
From the proof of Proposition \ref{Bryant}, we see that  
$(c(t)^2 - |\omega(t)|^2)\varphi\in\Omega^3_1$, 
$2c(t)\ast(\omega(t)\wedge\varphi)\in \Omega^3_{7}$, and 
$2\omega(t)\wedge\ast(\omega(t)\wedge\ast\varphi)
\in\Omega^3_1 \oplus \Omega^3_{27}$ with respect to 
$\varphi$, for all $t$.  Therefore derivatives will preserve the same irreducible subspaces induced by 
$\varphi$.  Applying this fact, the $\Omega^3_7$ - term in (\ref{metaji}) must be zero, which means that 
$\dot{c}(0)\,\omega_0 + c_0\,\dot{\omega}(0)=0$ by the definition of $\Omega^3_7$.  

\quad Consider points $(c_0,\omega_0)$ where
$c_0\ne 0$.  It follows that 
\begin{equation}\label{mini}
 \dot{\omega}(0)= -\frac{\dot{c}(0)}{c_0}\,\omega_0.
\end{equation} 
Note that this means if $\dot{c}(0)=0$ then 
$\dot{\omega}(0)=0$ as well, and we would be done with this case.  Thus we assume 
$\dot{c}(0)\ne 0$.  Differentiating 
$c(t)^2 + |\omega(t)|^2 =1$ and using (\ref{mini}) gives us 
\begin{align}
0 &= c_0\,\dot{c}(0)+\langle \dot{\omega}(0), \omega_0\rangle \notag\\
  &= c_0\,\dot{c}(0)-\frac{\dot{c}(0)}{c_0}\,|\omega_0|^2, \notag
\end{align}
or $c_0^{\,2} = |\omega_0|^2$.  On the other hand, from (\ref{mini}) and the proof of Proposition \ref{Bryant} we also have 
\begin{align}
 \dot{\omega}(0)\wedge\ast(\omega_0\wedge\ast\varphi) 
 \, +\, \omega_0\wedge\ast(\dot{\omega}(0)\wedge\ast\varphi) &= -2\frac{\dot{c}(0)}{c_0}\,
 \omega_0\wedge\ast(\omega_0\wedge\ast\varphi) \notag\\
 &= -2\frac{\dot{c}(0)}{c_0}\,
 \Big(\frac{3}{7}|\omega_0|^2\varphi + 
 \pi_{27}(\omega_0\wedge(\omega_0\wedge\ast\varphi))\Big).
\end{align}
Since the $\Omega^3_1$-term in (\ref{metaji}) must also be zero, it follows that 
\[
 c_0\,\dot{c}(0) - \frac{3\dot{c}(0)}{7c_0}
 |\omega_0|^2 = 0,
\]
or $c_0^{\,2} = 3|\omega_0|^2 / 7$.  We now get a contradiction to the previous conclusion of $c_0^{\,2} = |\omega_0|^2$ unless 
$\omega_0 =0$.  However, $\omega_0 = 0$ would imply $c_0 = 0$.  Therefore at $(c_0, \omega_0)$ where $c_0 \ne 0$, (\ref{metaji}) implies 
$(\dot{c}(0), \dot{\omega}(0)) = (0, 0)$ necessarily.

\quad Next we consider the case $c_0 = 0$. Then it follows by differentiating 
$c(t)^2 + |\omega(t)|^2 =1$  that 
$\langle \dot{\omega}(0), \omega_0\rangle=0$.  
Then by the standard identity in Proposition A.3 of \cite{K1} again, we see that
\begin{align}
 \dot{\omega}(0)\wedge\ast(\omega_0\wedge\ast\varphi)\wedge\ast\varphi 
 &= \dot{\omega}(0)\wedge\ast 3\omega_0 \notag\\
 &= 3\langle\dot{\omega}(0), \omega_0\rangle\,dv_\varphi
 \notag\\
 &=0 \notag.
\end{align}
Therefore $\pi_1(\dot{\omega}(0)\wedge\ast(\omega_0\wedge\ast\varphi))=0$, and hence 
$\dot{\omega}(0)\wedge\ast(\omega_0\wedge\ast\varphi)
\in\Omega^3_{27}$.  The identical argument shows that 
$\omega_0\wedge\ast(\dot{\omega}(0)\wedge\ast\varphi)
\in\Omega^3_{27}$ as well.  Now assuming (\ref{metaji}), it follows that 
\begin{equation}\label{question}
\dot{\omega}(0)\wedge
\ast(\omega_0\wedge\ast\varphi) 
 \, +\, \omega_0\wedge\ast(\dot{\omega}(0)
 \wedge\ast\varphi) = 0.
\end{equation} 
Recall that (\ref{metaji}) also implies $\dot{c}(0)\,\omega_0 + c_0\,\dot{\omega}(0)=0$, and so our assumption of $c_0 =0$ immediately implies $\dot{c}(0)=0$.  It remains to show that $\dot{\omega}(0)=0$ follows from (\ref{question}).  To this end, we first observe that 
\begin{equation}\label{referee}
\dot{\omega}(0)\wedge
\ast(\omega_0\wedge\ast\varphi) 
 \, +\, \omega_0\wedge\ast(\dot{\omega}(0)
 \wedge\ast\varphi) = h\odot\varphi,
\end{equation}
 where 
 $h = \dot{\omega}(0)\otimes\omega_0 + \omega_0\otimes\dot{\omega}(0)$.  To verify 
 (\ref{referee}), let $\{e_i\}_{i=1}^7$ be a point-wise frame adapted to the $G_2$-structure and 
 $\{\omega_i\}_{i=1}^7$ its metric dual frame.  Then write $\dot{\omega}(0)=\alpha_i\omega_i$ and 
 $\omega_0 = \beta_i \omega_i$, from which it follows that $h_{ij} = \alpha_i\beta_j + \beta_i\alpha_i$.  Then (\ref{referee}) follows directly from (\ref{Jazz}) in Section \ref{prelim}.  Now, by the isomorphism of $\text{Sym}^2_0(TM)$ onto 
 $\Omega^3_{27}$ via the map $h\mapsto h\odot\varphi$, (\ref{question}) implies $h=0$.  Then we compute 
 \begin{align}
  h(\dot{\omega}(0)^{\#},e_j) &= h_{ij}\alpha_i \notag\\
  &= (\alpha_i\beta_j + \beta_i\alpha_i)\,\alpha_i \notag \\
  &= |\dot{\omega}(0)|^2\beta_j + \langle \omega_0, \dot{\omega}(0) \rangle\,\alpha_j = |\dot{\omega}(0)|^2\beta_j.\notag
 \end{align}  
Since $\omega_0 \ne 0$, there must be some $\beta_j\ne 0$, which implies that $\dot{\omega}(0) =0$ necessarily. 
\qed



\quad      
Lastly, we want to reiterate that the diffeomorphsim 
$\Gamma\cong \mathbb{RP}^{\,b^1}$ is independent of 
the choice of $\{\varphi\}$ on $\mathcal{M}$, since it only depends on the first Betti number $b^1$ of $M$.  
Thus only the topology of the underlying manifold $M$ is relevant to the topology of $\Gamma$.

\section{The parameter space $N_f /G_2$.}\label{Lie}
For those familiar with the theory and language of principal bundles, it will be self-evident that most of the discussion in this section 
apply to all other holonomy groups.  Let $f$ be a frame over $p\in M$ that is 
adapted to some $G_2$-structure $\varphi$ - here thought of as a principal sub-bundle of the $GL(7,\mathbb{R})$-frame bundle over $M$ with structure group $G_2$.  It is well-known that $\varphi$ is torsion-free if and only if 
$\text{Hol}_f(D) \subset G_2$, where 
$D$ is the Levi-Civita connection associated to $g_\varphi$.  An equivalent characterization of $\varphi$ being torsion-free is that the $G_2$-structure is preserved by parallel-translation using $D$.  Therefore assuming 
$\text{Hol}_f(D) \subset G_2$, we can parallel translate $f$ to reconstruct the corresponding $G_2$-structure due to $M$ being connected.  In particular, a torsion-free $G_2$-structure (over $M$ connected) is determined by the fiber of the associated principal sub-bundle over a single point $p\in M$.  We denote this fiber by $[f]$, which is the orbit $\{f.h \,| \, h\in G_2\}$ with respect to the free and transitive right-action of the general linear group on the set of all (tangent) frames over a point.  Then note that 
two distinct $G_2$-structures over $M$ that generate the same Riemannian metric share the same Levi-Civita connection, which must preserve different frames due to the uniqueness property of parallel-translation.  Therefore the two $G_2$-structures must have distinct fibers over every 
point on $M$.


\quad Any frame over $p\in M$ adapted to the corresponding 
unique $SO(7)$-structure (the metric $g_\varphi$) can be written as $f.h$ for some unique $h\in SO(7)$.  Then for the frame $f$ above, by the well-known formula 
$\text{Hol}_{f.h} (D) = h^{-1}\text{Hol}_f(D) h$ it follows that 
$[f.h]$ represents a torsion-free 
$G_2$-structure also generating $g_\varphi$ if and only if $h^{-1}\text{Hol}_f(D) h \subset G_2$.  Therefore  
$[f.h]\in \Gamma$ if and only if $h^{-1}\text{Hol}_f(D)h \subset G_2$, and hence definition (\ref{N_f}) makes sense.  

\quad $N_f$ is a parameterization, based at $f$, of the set of all frames $\tilde{f}$ adapted to $g_\varphi$ (treated as an $SO(7)$-structure) 
such that $\text{Hol}_{\tilde{f}} (D) \subset G_2$.  
Suppose we had chosen a different frame $\tilde{f}$ adapted to another torsion-free $G_2$-structure that also generates $g_\varphi$, 
then by definition $\tilde{f} = f.g$ for some $g\in N_f$ and we see that 
\begin{align}
   h\in N_{\tilde{f}} \hskip 0.3cm &\Longleftrightarrow \hskip 0.3cm  
   (gh)^{-1}\text{Hol}_f(D) gh = h^{-1}\text{Hol}_{\tilde{f}}(D)h \subset G_2. \notag \\
   &\Longleftrightarrow gh \in N_f.
\end{align}   
The map $N_f \longrightarrow N_{\tilde{f}}$ given by $h\mapsto g^{-1}h$ is clearly bijective, 
and it furnishes a change of parameterization of the same set of frames via ``shifting'' by the element $g$.

\quad We define the map $N_f \longrightarrow \Gamma$ 
by $h\mapsto [f.h]$, and it is clearly surjective by definition.  Next we note that 
\begin{align}
[f.h] &= [f.g] \,\,\text{for} \,\, h,g\in N_f \notag\\  
      & \Longleftrightarrow f.h = (f.g).\tilde{h}=f.g\tilde{h} \,\,\text{for some}\, \tilde{h}\in G_2 \notag \\ 
      & \Longleftrightarrow g^{-1}h = \tilde{h}\in G_2 \notag\\
      & \Longleftrightarrow h\,G_2 = g\,G_2. \notag  
\end{align}
The second equivalence above follows from the fact that the action of the general linear group on ordered basis is free.  It then  follows that the induced map $N_f / G_2 \longrightarrow \Gamma$ given by $h\,G_2 \mapsto [f.h]$ is well-defined and bijective.  
Proposition \ref{prop1} is now verified.  We want to reiterate that 
in general $N_f / G_2$ is only a coset space, since $N_f$ may not be a group.  There is also the induced map 
$N_f / G_2 \longrightarrow N_{f.g} / G_2$, where the 
parametrization of $\Gamma$ by $N_{f.g}/G_2$ is equivalent to that by $N_f/G_2$ up to a shifting of cosets via 
left-multiplication by $g^{-1}$, for any $g\in G_2$.

\quad Note that according to (\ref{N_f}), 
when $\text{Hol}_f(D) = G_2$, $N_f$ is exactly the normalizer $N(G_2)$ of $G_2$ in $SO(7)$.  Since normalizers are 
Lie groups, $N(G_2)/G_2$ is also a Lie group.  
In \cite{Kz}, it was shown by purely algebraic means that $N(G_2) = G_2$.  Therefore by Proposition \ref{prop1} we have the following result. 
\newtheorem{corollary}{Corollary}[section]
 \begin{corollary}
  Let $M$ be a connected $7$-dimensional manifold admitting a $G_2$-structure $\varphi$ with full $G_2$-holonomy.  
  Then $\varphi$ is the only torsion-free $G_2$-structure 
  on $M$ that generates the metric $g_\varphi$.   
 \end{corollary}
This is the same result that was obtained in the previous section from Bryant's formula (\ref{claim}).  It is also interesting here 
that a result in pure algebra implies a result in moduli spaces by way of the geometric structures involved.  

\quad On the other hand, 
we can also 
turn the argument around and give a new moduli space argument of why $N(G_2) = G_2$.  
Set-theoretically, we have $N_f / G_2 \cong \Gamma \cong \mathbb{RP}^{\,b_1}$ by incorporating the result from Section \ref{Bry}.  By taking 
any one of the compact $G_2$-manifolds with full $G_2$-holonomy constructed by Joyce (see \cite{J1}), we have $b_1=0$ and $N_f = N(G_2)$.  Then 
$N(G_2) / G_2 \cong \mathbb{RP}^{\,0} = \{\text{single point}\}$, which implies that $N(G_2) = G_2$ necessarily.

\quad It can be shown that $N_f$ is a submanifold of $SO(7)$.  Therefore $N_f / G_2$, which 
is the orbit space of $G_2$ acting on $N_f$ by right-multiplication, is also a manifold.  We now argue that 
$N_f / G_2$ is diffeomorphic, and hence homeomorphic to $\Gamma$. 

\newtheorem{prop4}[corollary]{Proposition}
\begin{prop4}
 The map $N_f/G_2 \longrightarrow \Gamma$ given by $hG_2\mapsto [f.h]$ is a diffeomorphism.  
\end{prop4}

\text{\bf{Proof.}}\quad We already saw that 
$hG_2\mapsto [f.h]$ is a bijection.  We will show that this map is in fact locally a diffeomorphism, via the inverse function theorem.  

\quad Consider $\tilde{f} = f.g$ for any $g\in N_f$.  
Let $h(t)$ be a smooth family of elements in $N_{\tilde{f}}$ such that $h(0) = I$.  Up to parallel translation on $M$, we can write 
$[\tilde{f}.h(t)]=\big(\tilde{f}.h(t)\big)^*\varphi_0 = 
\varphi(t)$ as 
elements in $\Gamma$.  We also denote $\tilde{\varphi}= 
\tilde{f}^*\varphi_0$.  By the natural action of the general linear group on frames and the fact that 
$h(t)\in SO(7)$, it follows that 
$\tilde{f}.h(t) = h(t)^{-1}\tilde{f}$, where the right-hand side denotes a composition of linear maps.  Then for 
tangent vectors $u, v, w$ we see that 
\begin{align}\label{showered}
 \varphi(t) (u,v,w) &= \big(\tilde{f}.h(t)\big)^*\varphi_0
 (u,v,w)\notag\\
 &=\varphi_0 \big(h(t)^{-1}\tilde{f}(u), 
 h(t)^{-1}\tilde{f}(v), 
 h(t)^{-1}\tilde{f}(w)\big) \notag\\
 &= \tilde{f}^*\varphi_0 \big(\tilde{f}^{-1}h(t)^{-1}
 \tilde{f}(u), \tilde{f}^{-1}h(t)^{-1}\tilde{f}(v), 
 \tilde{f}^{-1}h(t)^{-1}\tilde{f}(w)\big) \notag\\
 &= \tilde{\varphi}\big(\tilde{f}^{-1}h(t)^{-1}
 \tilde{f}(u), \tilde{f}^{-1}h(t)^{-1}\tilde{f}(v), 
 \tilde{f}^{-1}h(t)^{-1}\tilde{f}(w)\big).
\end{align} 
Then
\begin{align}\label{AHC}
 \frac{\partial \varphi(t)}{\partial t}(u,v,w)\Big|_{t=0} &= \tilde{\varphi}\big(-\tilde{f}^{-1}h'(0)\tilde{f}(u),v,w\big) + \tilde{\varphi}\big(u, -\tilde{f}^{-1}h'(0)\tilde{f}(v),w\big) \notag\\
 &\hskip 0.3cm +\tilde{\varphi}\big(u, v, -\tilde{f}^{-1}h'(0)\tilde{f}(w)\big)\notag\\
 &= \big(-\tilde{f}^{-1}h'(0)\tilde{f}\,\big) \odot \tilde{\varphi}
 (u,v,w),
\end{align} 
where $(h^{-1})^{'}(0)= -h'(0)$ follows from $h(0)=I$.

\quad Now, suppose $\frac{\partial \varphi(t)}
{\partial t}\big|_{t=0}=0$.  
By (\ref{AHC}) and the discussion following (\ref{split}) in 
Section \ref{prelim}, we conclude that $h'(0)$ must 
be an element of the Lie algebra of $G_2$.  This shows 
that for any transverse slice through $I\in G_2$ in $N_{\tilde{f}}$, the derivative of the map $h\mapsto [\tilde{f}.h] = \big(\tilde{f}.h\big)^*\varphi_0$ is injective.  By the inverse function theorem this shows that the map is locally invertible onto the transverse slice, with the inverse also differentiable.  Passing to cosets this shows that the induced map $N_{\tilde{f}}/G_2\longrightarrow \Gamma$ is locally invertible about the identity element $G_2$ with a differentiable inverse.

\quad Next we refocus on the map 
$N_f/G_2 \longrightarrow \Gamma$.  For any $g\in N_f$, 
we have $gG_2 \mapsto [f.g]=[\tilde{f}]=\tilde{f}^*\varphi_0 = \tilde{\varphi}$, using the notations above.  
Then note that the assignment above can be factored via the composition 
\[
 N_f/G_2 \longrightarrow N_{\tilde{f}}/G_2\longrightarrow \Gamma
\]  
of maps, where we recall that $N_f/G_2 \longrightarrow N_{\tilde{f}}/G_2$ is given by left-multiplication by $g^{-1}$ and is clearly a diffeomorphism.  Then 
since we saw that the map $N_{\tilde{f}}/G_2\longrightarrow \Gamma$ 
is a local diffeomorphism about $G_2$, this shows that $N_f/G_2 
\longrightarrow \Gamma$ is a diffeomorphism.    
\qed

\quad Note that a-priori, $N_f$ (and hence $N_f / G_2$) could have multiple 
connected components.  The work above however, shows that the manifold $N_f / G_2$ actually has exactly one component.  
Also, note that we identified the topology of $N_f / G_2$ by analyzing relevant geometric structures, and not through direct Lie group-theoretic means.

\section{Implications on the Moduli Space}\label{Moduli}
Consider the continuous projection map $\Gamma \longrightarrow \mathcal{M}$ given by $\tilde{\varphi} \mapsto \{\tilde{\varphi}\}$.  
If this map is injective, then the topology in the $H^3_7(M)$-directions on $\mathcal{M}$ is given by a space homeomorphic to $\mathbb{RP}^{\,b_1}$.  However, it is conceivable that for some $\tilde{\varphi}\in \Gamma$, 
$f^{*}\tilde{\varphi}$ is 
also in $\Gamma$ for some $f\in\text{Diff}_0(M)$ such that $f^{*}\tilde{\varphi}\ne\tilde{\varphi}$.  Note that in this case, $f$ is necessarily an isometry of $g_{\tilde{\varphi}}$ since it can be shown that $g_{f^{*}\tilde{\varphi}} = 
f^{*}g_{\tilde{\varphi}}$ for any diffeomorphism 
$f$ of $M$ and any $G_2$-structure 
$\tilde{\varphi}$.   In other words, 
$\Gamma$ could intersect a given orbit of $\mathcal{M}$ more than once.  Throughout this section $f$ denotes a diffeomorphism of $M$, and not a point-wise frame as in Section \ref{Lie}.

\quad The observation above points to the structure of covering spaces.   Let $X$ denote the image of the map $\Gamma \longrightarrow \mathcal{M}$ above.  The topology of $X$ is exactly what is meant by the topology of $\mathcal{M}$ in directions tangent to $H^3_7(M)$. Since $\tilde{\varphi}$ and 
$f^*\tilde{\varphi}$ projects onto the same point in $X$ for any isometry of $g_{\tilde{\varphi}}$ in $\text{Diff}_0(M)$, it is natural to expect that $X$ is the orbit space of the action by the pull-back of such isometries on $\Gamma$, and that $\Gamma\longrightarrow X$ is a covering map. 

\quad We begin with the following definition and the ensuing basic result.  Given a differential form $\eta$, we say that a smooth family $f_t$ of diffeomorphisms 
of $M$ is \text{\it{regular}} (relative to $\eta$) if 
$\frac{\partial}{\partial t}f_{t}^{*}\eta \ne 0$ for all $t$.   

\newtheorem{lem}{Lemma}[section]
\begin{lem}\label{iso}
 On a compact $7$-manifold with torsion-free $G_2$-structure 
 $\varphi$, there does not exist any 
 smooth, regular (relative to $\varphi$) family of isometries of $g_\varphi$.       
\end{lem}

\text{\bf{Proof.}}\quad Let $f_t$ be a smooth family 
of isometries of $g_\varphi$.  Then 
$f_t^{*}\varphi$ is a family of torsion-free $G_2$-structures on $M$.  From the discussion in Section \ref{Bry}, 
we know that the velocity vectors of the curve 
$f_t^{*}\varphi$ in $\Omega^3(M)$ are 
harmonic $\Omega^3_7$-forms.  On the other hand, we see 
that 
\begin{equation}\label{transv}
\frac{\partial f_t^{*}\varphi}{\partial t} 
= f_t^{*}L_V \varphi = 
L_{(f_t)_{*}V}\varphi = 
d((f_t)_{*}V\lrcorner\varphi),        
\end{equation}
which is exact and hence equals $0$ if harmonic.  This simple 
calculation implies that 
$f_t^{*}\varphi = f_0^{*}\varphi$ for all $t$.  
Thus $f_t^{*}$ cannot be a regular family relative to $\varphi$.     \qed   

\quad We denote by $S_0$ the set of all isometries of $g_\varphi$ in $\text{Diff}_0(M)$.  
Note that because $g_{f^{*}\tilde{\varphi}} = 
f^{*}g_{\tilde{\varphi}}$, we have  
$f^{*}:\Gamma\longrightarrow \Gamma$ for every 
$f\in\Sigma$.  Define the binary operation  
$f\#\tilde{f} = \tilde{f}\circ f$ for all $f,\tilde{f}\in S_0$.  The reversing of order of functions is tailored to giving a well-defined left group action of $S_0$ on $\Gamma$, as we shall see below.  

\begin{lem}\label{S_0}
 $S_0$ is a compact Lie group with respect to the operation $\#$.  Moreover, pull-back of elements in  $\Gamma$ by elements in $S_0$ is a left group action, so that $X$ is exactly the orbit space of this action.
\end{lem}

\text{\bf{Proof.}}\quad First we verify the group structure of $S_0$.  Clearly, the identity $I\in S_0$.  Next, $f\#\tilde{f} = \tilde{f}\circ f$ is still an isometry for any $f,\tilde{f}\in S_0$.  Furthermore, $f\#\tilde{f}\in \text{Diff}_0(M)$ because $\text{Diff}_0(M)$ is a group. Note that $f\#\tilde{f}$ is not necessarily isotopic to the identity through isometries - as only diffeomorphisms are required by the definition of $S_0$.  Thus $f\#\tilde{f}\in S_0$ as 
well.  Clearly the identity $I$ is in $S_0$, and that $I\# f = f\# I = f$ for all $f\in S_0$.  For any $f\in S_0$, its inverse $f^{-1}$ is also an isometry and simultaneously a member of $\text{Diff}_0(M)$, and so $S_0$ is also closed under taking inverses.  Lastly, the associativity of $\#$ is inherited from $\text{Diff}_0(M)$ as well.  Therefore $S_0$ is indeed a group.                     

\quad Next we establish the smooth structure on $S_0$.  Recall that the $7$-manifold $M$ is compact.  It is well-known that the isometry group $\text{Isom}(g)$ of any metric $g$ on a compact manifold $M$ is a compact Lie group.  Therefore $\text{Isom}(g_\varphi)$ is a compact Lie group.  Also note that since $\text{Diff}_0(M)$ is the connected component of $\text{Diff}(M)$ containing the identity, it follows from general topology that $\text{Diff}_0(M)$ must be a closed subset of $\text{Diff}(M)$.  On the other hand, $\text{Diff}(M)$ is also Hausdorff, therefore the compact subset $\text{Isom}(g_\varphi)$ must also be closed in $\text{Diff}(M)$ as well.  This means 
$S_0 = \text{Isom}(g_\varphi)\cap\text{Diff}_0(M)$ is also closed in the compact set $\text{Isom}(g_\varphi)$, hence $S_0$ is also compact.  Finally, since $S_0$ is a (topologically) closed subgroup of $\text{Isom}(g_{\varphi})$, it inherits the smooth structure of $\text{Isom}(g_\varphi)$.  This completes the proof that $S_0$ is a compact Lie group.     

\quad We define a left group action of $S_0$ on $\Gamma$ by $f.\tilde{\varphi} = f^*\tilde{\varphi}$.  Let us check that this is indeed a left group action.  Consider any $\tilde{\varphi}\in \Gamma$.  Clearly
$I.\tilde{\varphi} = \tilde{\varphi}$.  Then we see that 
\[
 \tilde{f}.(f.\tilde{\varphi}) = \tilde{f}.(f^*\tilde{\varphi}) 
 = \tilde{f}^* f^*\tilde{\varphi} 
 = (f\circ\tilde{f})^*\tilde{\varphi} = 
 (\tilde{f}\# f).\tilde{\varphi}
\]  
for any $f,\tilde{f}\in S_0$.  This verifies the group action.  The space $X$ is then the orbit space of this action by the definition of the projection $\Gamma\longrightarrow X$.  \qed

\newtheorem{cor5}[lem]{Corollary}
\begin{cor5}
Let $M$ be a compact $7$-manifold with torsion-free $G_2$-structure $\varphi$.  If the only isometry of $g_\varphi$ in $\text{Diff}_0(M)$ is the identity, then  $X$ 
 is homeomorphic to $\mathbb{RP}^{\,b^1}$.   
\end{cor5}

\quad To show that $\Gamma\longrightarrow X$ is a covering map, it would be sufficient to show that 
$S_0$ acts \text{\it{properly discontinuously}} on $\Gamma$, which would then imply that $S_0$ is the group of deck transformations.  However, when $b_1 \geq 1$ there exist \text{\it{symmetries}} of the underlying $G_2$-structure (see \cite{Lin}), which means the action by $S_0$ is not free and hence cannot be properly discontinuous.  We will elaborate on this in more detail later, but in spite of it we can still show that $\Gamma\longrightarrow X$ is a covering map.  To do so, we will examine the nature of each orbit of the action by $S_0$.  In what follows, we will denote by ``$\text{orb}_{S_0}(\tilde{\varphi})$'' as the orbit containing $\tilde{\varphi}$ of $S_0$ acting on $\Gamma$, and denote by ``$\text{stab}_{S_0}(\tilde{\varphi})$'' as the stabilizer subgroup of $S_0$ fixing $\tilde{\varphi}$.       

\begin{lem}\label{finite}
 For every $\tilde{\varphi}\in\Gamma$, $\text{orb}_{S_0}(\tilde{\varphi})$ is a finite subset of $\Gamma$ and $\text{stab}_{S_0}(\tilde{\varphi})$ is a compact Lie subgroup of $S_0$ with dimension $b^1$.  In particular, $S_0$ is also of dimension $b^1$. 
\end{lem}

\text{\bf{Proof.}}\quad 
Since $S_0$ is compact by Lemma \ref{S_0}, it is also a \text{\it{proper}} action.  Then by the standard fact that orbits for a proper action on any manifold are closed submanifolds, it follows that $\text{orb}_{S_0}(\tilde{\varphi})$ is a closed submanifold of $\Gamma$ for any $\tilde{\varphi}\in\Gamma$.  Then by Lemma \ref{iso}, $\text{orb}_{S_0}(\tilde{\varphi})$ must be a $0$-dimensional submanifold of $\Gamma$.  In particular, $\text{orb}_{S_0}(\tilde{\varphi})$ is just a finite collection of points in $\Gamma$, otherwise (via the compactness of $\Gamma$) a limit point would violate the definition of $\text{orb}_{S_0}(\tilde{\varphi})$ being a submanifold of $\Gamma$.   

\quad The stabilizers of proper actions are compact, so it follows that $\text{stab}_{S_0}(\tilde{\varphi})$ is compact.  In particular, since $\text{stab}_{S_0}(\tilde{\varphi})$ is also a closed subset of of the Lie group $S_0$, it must also be a Lie subgroup of $S_0$.  Recall that $\text{orb}_{S_0}(\tilde{\varphi}) \cong S_0/\text{stab}_{S_0}(\tilde{\varphi})$ as sets.  Thus the cosets of $\text{stab}_{S_0}(\tilde{\varphi})$ must be a finite collection of mutually-disjoint compact submanifolds of $\text{Isom}(g_\varphi)$, and each coset is a translation of (hence homeomorphic to) $\text{stab}_{S_0}(\tilde{\varphi})$.  In \cite{Lin} it was shown that $b^1$ is the dimension of the vector space of (continous) symmetries of $\varphi$. Therefore by the standard translation argument we see that as a manifold, $\text{stab}_{S_0}(\tilde{\varphi})$ (hence each of its cosets) must also have dimension equal to $b^1$.  Finally, since $S_0$ is also the disjoint union of the cosets of 
$\text{stab}_{S_0}(\tilde{\varphi})$, it follows that $S_0$ must also have the same dimension $b^1$.       \qed

\text{\bf{Proof of Theorem \ref{Thm1}.}}\quad 
Let $\{\tilde{\varphi}\}\in X$.  Then by definition $\text{orb}_{S_0}(\tilde{\varphi})$ is the pre-image of the projection $\Gamma\longrightarrow X$. The finiteness of $\text{orb}_{S_0}(\tilde{\varphi})$ from Lemma \ref{finite} allows us to choose an open neighborhood $U$ in $X$ containing $\{\tilde{\varphi}\}$, so that there are mutually-disjoint open neighborhoods in $\Gamma$ containing each point in $\text{orb}_{S_0}(\tilde{\varphi})$ which project onto $U$.  This concludes the proof. \qed  

Note that by the smoothness of $\mathcal{M}$, the number of sheets in 
the covering map $\Gamma\longrightarrow X$ 
is constant over each connected component of $\mathcal{M}$, but may vary across different components. 

\quad  When $b^1 = 0$, both $\Gamma$ and $X$ are respectively just the single point $\varphi$.  Then $\Gamma \longrightarrow X$ is the trivial covering map with the identity as its only deck transformation.  Since $\text{orb}_{S_0}(\varphi) \cong S_0/\text{stab}_{S_0}(\varphi)$ as sets, this implies that every element of $S_0$ fixes a $G_2$-structure $\varphi$ with full $G_2$-holonomy.   Moreover, $S_0$ is a $0$-dimensional submanifold of the compact Lie group $\text{Isom}(M)$ by Lemma \ref{finite}, so it must be finite.  Next, we investigate the situation where $b^1 = 1$.  From this point on, we also denote by $\Sigma$ the group of deck transformations of $\Gamma\longrightarrow X$.  

\newtheorem{prop5}[lem]{Proposition}
\begin{prop5}\label{38}
If $b^1 = 1$, then $X$ is homeomorphic to $S^1$ for every 
$\varphi$ that represents a point on $\mathcal{M}$.  Moreover, $\Sigma\cong\mathbb{Z}_q$ for some positive integer $q$.  
If in addition $b^3_{27}=0$, then each connected component of $\mathcal{M}$ 
is homeomorphic to the punctured-plane.    
\end{prop5}  

\text{\bf{Proof.}} \quad Since $\mathbb{RP}^1$ is 
homeomorphic to $S^1$, the corresponding covering map 
can be written as $S^1 \longrightarrow X$.  The covering map 
is continuous, so $X$ must be compact, connected, and $1$-dimensional (locally as a differentiable manifold).  Because 
the local deformation of $\Gamma$ is tangent to the 
$1$-dimensional subspaces $H^3_7(M)$, there cannot be any self-intersections of $X$ as a submanifold in 
$\mathcal{M}$.  Thus $X$ must be homeomorphic to $S^1$ for every $\varphi$ that represents a point on 
$\mathcal{M}$.  It follows that if in addition $b^3_{27}=0$, then $\mathcal{M}$ is a 
$2$-dimensional manifold necessarily homeomorphic to 
$(0,\infty)\times S^1$, or the punctured-plane.

\quad  Let $H$ denote the (injective) push-forward  
of $\pi_1(\Gamma)$ under the 
covering map $\Gamma\longrightarrow X$. Recall from the theory of covering spaces that the group of deck transformations is isomorphic to 
the normalizer of $H$ in $\pi_1(X)$ modulo $H$.  
Since in our case $\pi_1(X) \cong \pi_1(S^1) \cong \mathbb{Z}$, 
hence ablelian, the group of deck transformations is isomorphic to $\pi_1(X) / H$.  In other words, the covering map $\Gamma\longrightarrow X$ is \text{\it{normal}} (see \cite{H}).  On the other hand, $\pi_1(\Gamma)\cong\pi_1(S^1)\cong\mathbb{Z}$ as well, so $H$ must also be isomorphic to a nontrivial subgroup of $\pi_1(X)\cong\mathbb{Z}$.  Elementary group theory tells us that such a subgroup must be $q\mathbb{Z}$ for some positive integer $q$.  Taking the quotient, it follows that 
the group of deck transformations is isomorphic to $\mathbb{Z} / q\mathbb{Z}\cong\mathbb{Z}_q$.   \qed    

  \quad Note that a-priori $\text{stab}_{S_0}(\varphi)$, and hence $S_0$, could be disconnected.  Based on Lemma \ref{S_0} and Lemma \ref{finite}, we expect $S_0$ to be dijoint copies of circles.  Therefore we make the following conjecture:

\text{\bf{Conjecture.}}\quad If $b^1 = 1$, then $S_0$ is isomorphic to $S^1\times \mathbb{Z}_q$ for some positive integer $q$. 

\quad For $b^1 > 1$, it is natural to consider the composition $S^{\,b^1}\longrightarrow\mathbb{RP}^{\,b^1}\longrightarrow X$, with the standard covering map $S^{\,b^1}\longrightarrow\mathbb{RP}^{\,b^1}\cong \Gamma$ that has $\mathbb{Z}_2$ as its group of deck transformations.  Since $S^{b^1}$ is simply-connected, the resulting universal covering $S^{\,b^1}\longrightarrow X$ is normal, and its group $G$ of deck transformations is isomorphic to $\pi_1(X)$.  Moreover, $X$ is homeomorphic to 
the orbit space $S^{\,b^1}/G$, where $G$ acts freely on $S^{\,b^1}$.  

\quad Note that if $\Gamma \longrightarrow X$ is also a normal covering, then $\Sigma$ is isomorphic to $\pi_1(x)$ modulo the injective push-forward of $\pi_1(\Gamma) \cong \pi_1(\mathbb{RP}^{\, b^1})\cong \mathbb{Z}_2$.  It then follows that $G$ is a group extension of $\Sigma$ by $\mathbb{Z}_2$, i.e. there exists a short exact sequence 
\[
 0 \longrightarrow \mathbb{Z}_2 \longrightarrow G 
 \longrightarrow \Sigma \longrightarrow 0.
\] 
If $\Gamma$ is a normal covering space over $X$, the image of 
the $2$-element group $\mathbb{Z}_2$ must be contained in the center of $G$.  Therefore $G$ is a \text{\it{central extension}} - an algebraic object classified by the group cohomology $H^2(\Sigma, \mathbb{Z}_2)$.  Note that if $G$ splits, it can only do so as the trivial extension $\mathbb{Z}_2 \times \Sigma$\, since the automorphism group of $\mathbb{Z}_2$ is trivial, and this corresponds to the zero element in $H^2(\Sigma, \mathbb{Z}_2)$.
  
\quad Next we observe that in the proof of 
Lemma \ref{finite}, $\text{orb}_{S_0}(\tilde{\varphi})$ is actually diffeomorphic to the quotient manifold $S_0/\text{stab}_{S_0}(\tilde{\varphi})$.  Therefore if we assume $\text{stab}_{S_0}(\tilde{\varphi})$ is a normal subgroup of $S_0$, then $S_0/\text{stab}_{S_0}(\tilde{\varphi})$ is a Lie group that acts transitively on $\text{orb}_{S_0}(\tilde{\varphi})$.  It would then follow that $X$ is diffeomorphic to the orbit space of $S_0/\text{stab}_{S_0}(\tilde{\varphi})$ acting on $\Gamma$, which also means that $\Gamma \longrightarrow X$ would be a normal covering. In particular, it follows that 
$S_0/\text{stab}_{S_0}(\tilde{\varphi})$ is exactly the group $\Sigma$ of deck transformations of $\Gamma\longrightarrow X$.  Although it certainly seems artificial at the present, we shall conveniently assume that $\text{stab}_{S_0}(\tilde{\varphi})$ is a normal subgroup of $S_0$ for the cases $b^1 > 1$. 
 
\quad There are some partial results that can be obtained regarding the group extension problem above.  For example, by the Schur-Zassenhaus Lemma (see \cite{Rot} for reference), 
if $\Sigma$ is finite and whose order $|\Sigma|$ is relatively prime with $|\mathbb{Z}_2|=2$, then $H^2(\Sigma, \mathbb{Z}_2) = 0$.  
Then since $\Sigma$ has to be a finite group due to Lemma \ref{finite}, we can draw the following conclusion.
\begin{prop5}\label{ny}
 Suppose $b^1 >1$, $\text{stab}_{S_0}(\tilde{\varphi})$ is normal in $S_0$, and $\Sigma$ has odd order.  Then $X$ is homeomorphic to the orbit space of an action of 
 $\mathbb{Z}_2 \times \Sigma$ on $S^{\,b^1}$.
\end{prop5}  

\quad  One can also approach the problem of describing $X$ by asking which groups can act freely on spheres.  For finite groups it is known that such a group must have \text{\it{periodic cohomology}} - a condition equivalent to all abelian subgroups of the group being cyclic.  Via the K\"{u}nneth formula it can be shown that $\mathbb{Z}_p \times \mathbb{Z}_p$ does not have periodic cohomology for any prime $p$.  For a survey of these results we refer to 
\cite{Wall}. 

\begin{prop5}\label{la}
 Suppose $b^1 >1$, $\text{stab}_{S_0}(\tilde{\varphi})$ is normal in $S_0$, and $\Sigma \cong \mathbb{Z}_2$.  Then 
 $X$ is homeomorphic to the orbit of an action of $\mathbb{Z}_4$ on $S^{\,b^1}$.
\end{prop5}
\text{\bf{Proof.}}\quad $\Sigma \cong \mathbb{Z}_2$ implies that $G$ must be of order $4$, since 
$G / \mathbb{Z}_2 \cong \Sigma$.  The only possibilities are then $G\cong \mathbb{Z}_2\times\mathbb{Z}_2$ or $\mathbb{Z}_4$.  The former can be ruled out because not all its subgroups are cyclic.  Therefore the proposition follows.
\qed

On the other hand, the author does not know whether or not $\mathbb{Z}_4$ can act freely on spheres of odd dimension at least $3$.  The assumptions in Propositions \ref{ny} and \ref{la} also seem artificial at the present.  In particular, 
we would like there to be a way to determine whether or not $\text{stab}_{S_0}(\tilde{\varphi})$ is normal in $S_0$.  It would also be very helpful if there is some way of determining 
$\Sigma$ independently.  Ultimately, we should be able to understand the topology of $X$ more if our knowledge about free actions on spheres or the structure of $H^2(\Sigma,\mathbb{Z}_2)$ advances.         

\quad The remaining mystery concerning the topology of 
$\mathcal{M}$ lies in directions tangent to 
$H^3_{27}(M)$.  If $b_{27}^3\ne 0$, there may be some associated global topology like that of $\Gamma$ described in this paper.  On the other hand, according to (\ref{coho}) we have $b^3_{27} = b^3 - 1-b^1$.  In view of Proposition \ref{38}, the moduli space $\mathcal{M}$ of the example of $S^{1}\times Y$ in Section \ref{prelim} contains families of circles.  Note that by the K\"{u}nneth formula, we have 
$H^3(S^1\times Y) = H^2(Y)\oplus H^3(Y)$.  We may be able to find a Calabi-Yau $3$-fold such that $b^2(Y)+b^3(Y)= 2$, so that the moduli space $\mathcal{M}$ of 
$S^1\times Y$ is homeomorphic to the punctured plane.  This, and the many questions raised above, serve as a good source for future research.

\vskip 1.5cm

\text{\sc{Department of Mathematics and Statistics}} \\ 
\text{\sc{University of South Alabama}} \\
\text{\sc{Mobile, AL 36688}} \\
\text{\it{email}}: cclin@southalabama.edu

\end{document}